\newcount\notenumber

\def\note{\advance\notenumber by 1
\footnote{$^{(\the\notenumber)}$}}

\def\Pr{{\bf P}}

\def\Q{{\bf Q}}

\def\Z{{\bf Z}}

\def\x{{\bf x}}

\def\barQ{{\bar{\bf Q}}}
\def\Gm{{{\bf G}_m^2}}

\def\Z{{\bf Z}}

\def\OS{{\cal O}_S}

\def\barQ{{\overline{\bf Q}}}

\font\title=cmr10 scaled 1200

\centerline{{\title A lower bound for the height of a
rational function at $S$-unit points}}
\bigskip

\centerline{Pietro Corvaja\hskip 2truecm Umberto Zannier}
\bigskip

\noindent{\bf Abstract}. Let $a,b$ be given multiplicatively independent positive
integers and let $\epsilon >0$. In a recent paper written jointly also with Y. Bugeaud we
proved  the  upper bound 
$\exp(\epsilon n)$ for $\gcd (a^n-1,b^n-1)$; shortly afterwards we generalized this to the
estimate $\gcd (u-1,v-1)<\max(|u|,|v|)^\epsilon$, for multiplicatively independent $S$-units
$u,v\in\Z$.  In a subsequent analysis of those results it turned out that a perhaps  better 
formulation of them may be obtained in terms of the language of heights of algebraic
numbers. In fact, the purposes of the present paper are:
to generalize the upper bound for the g.c.d. to pairs of rational functions
other than $\{u-1,v-1\}$  and 
to extend the results to the realm of algebraic
numbers, giving at the same time a new formulation of the bounds in terms of height
functions and algebraic subgroups of $\Gm$. \bigskip

AMS Classification: 11J25.

Key words: Lower bounds for the height, Subspace Theorem, Linear tori.

\bigskip

\noindent{\bf \S 1 Introduction}. Let $a,b$ be given multiplicatively independent positive
integers and let $\epsilon >0$. In the paper [1] a nearly best-possible upper bound 
$\exp(\epsilon n)$ for g.c.d.$(a^n-1,b^n-1)$ was proved, 
as the integer $n$ tends to infinity. In
the subsequent paper [6] we observed (see Remark 1) that the combination of those 
arguments with [1]   yields a generalization, where $a^n, b^n$ are replaced by
multiplicatively independent $S$-units $u,v\in \Z$. Namely the upper bound
g.c.d.$(u-1,v-1)<\max(|u|,|v|)^\epsilon$ holds with finitely many exceptions. (Full details 
for this can be obtained in a straightforward way from [1] and [6], as is done in [13] and
[4].)

In a subsequent analysis of those results it turned out that a perhaps  better 
formulation of them may be obtained in terms of the language of heights of algebraic
numbers. This also suggested an extension  to general  number fields. Finally, in
view of possible applications (like e.g. the one in [6], where a conjecture about
the greatest prime factor of $(ab+1)(ac+1)$ was settled) we have tried to replace
$u-1,v-1$ with more general functions. 

In fact, the purposes of the present paper are:

 (i) {\it to generalize the upper bound for the g.c.d. to pairs of rational functions
other than $\{u-1,v-1\}$} and 

(ii) {\it to extend the results to the realm of algebraic
numbers, giving at the same time a new formulation of the bounds in terms of height
functions and algebraic subgroups of $\Gm$}.  \medskip

More precisely, note that a ``good" upper bound for the g.c.d. of two integers amounts
to a lower bound for  the height of their ratio, essentially by their maximum height.
With this in mind, we seek lower bounds for the height $h(\varphi (u,v))$, where
$\varphi$ is a rational function and $u,v$ run along   $S$-units in a number
field, for a finite set $S$. Equivalently, we may let $(u,v)$ run through a finitely
generated group  $\Gamma\subset\Gm(\barQ)$.  Our main results will give, under
suitable assumptions,  the expected lower bounds, with the possible exceptions of
the pairs $(u,v)$ contained in a finite union of proper algebraic subgroups (or
translates) of $\Gm$: these kind of subvarieties are certainly relevant in the
context; for instance, a pair $(u,v)$ lies in some proper algebraic subgroup of $\Gm$
if and only if $u,v$ are multiplicatively dependent. 

Before stating the results we observe that analogous investigations, from the
quantitative point of view appear in [2] (for the gratest prime factor of
$(ab+1)(ac+1)(bc+1)$) and [3] (where some results of [7] are in part quantified).\medskip

We immediately proceed to give the formal statements.

\bigskip

\noindent{\bf Theorem 1}. {\it Let $\Gamma\subset \Gm(\barQ)$
be a finitely generated group, $p(X,Y),q(X,Y)\in
\barQ [X,Y]$ be non constant coprime polynomials and suppose that
not both of them vanish at $(0,0)$. 
For every positive $\epsilon$, the Zariski closure of the set of
solutions  $(u,v)\in\Gamma$ to the inequality
$$
h\left ({p(u,v)\over q(u,v)}\right)<
h(p(u,v):q(u,v):1)- \epsilon\cdot \max\{h(u),h(v)\}\eqno(1.1)
$$
is the union of  finitely many translates of 1-dimensional subtori of $\Gm$ which can be
effectively determined and a finite set. }
\bigskip

Note that, since $h(x:y:1)\ge \max \{h(x),h(y)\}$, in (1.1) we can replace
$h(p(u,v): q(u,v):1)$ with $\max \{h(p(u,v)),h(q(u,v))\}$ (but this would lead to a weaker
statement). Also,  we could have replaced (1.1) with  the equivalent
inequality  
$$ 
\sum_{\mu\in M_K}\log^{-}\max\{|p(u,v)|_\mu,|q(u,v)|_\mu\}
<-\epsilon \max\{h(u),h(v)\},\eqno(1.2)
$$
where $\log^-(\cdot)$ stands for $\min\{0,\log(\cdot)\}$. (See the beginning of next
section for conventions about valuations.)

Actually, the proof will go through (1.2) and we shall show later the equivalence.
Note that  the left side of (1.2) is an analogue of the g.c.d. (in a sense which considers
also archimedean valuations), and for rational integers it is precisely its logarithm.
\medskip

Theorem 1 easily implies the following   corollary:
\medskip

\noindent {\bf Corollary 1.} {\it Let $\Gamma$ be as before.
Then for multiplicatively independent pairs $(u,v)\in \Gamma$
the height of the ratio $(u-1)/(v-1)$ verifies the asymptotic equivalence, 
as $\max \{h(u),h(v)\}\rightarrow\infty$:
$$
h((u-1)/(v-1))\sim h(1:u:v).
$$
}
\medskip

Note that, as before, this  is (essentially) equivalent to the upper bound
$$
\sum_{\mu\in M_K}\log^{-}\max\{|u-1|_\mu,|v-1|_\mu\}>- \epsilon\cdot
\max\{h(u),h(v)\},\eqno(1.3)
$$
valid for all $\epsilon$ and all multiplicatively independent pairs 
$(u,v)\in\Gamma$ apart a finite set depending on $\epsilon$ and $\Gamma$.

Also,  the main results of [1] and [6] 
are immediate consequences  of this Corollary.
\bigskip


Our most general result is the following
\medskip

\noindent{\bf Main Theorem}. {\it Let $f(X,Y)\in \barQ(X,Y)$ 
be a rational function
and $\Gamma\subset \Gm(\barQ)$ be a finitely generated subgroup.
Denote by $T_1(X,Y),\ldots,T_N(X,Y)$ the monomials appearing in 
the numerator and denominator of $f(X,Y)$ and suppose that 
$1\in\{T_1,\ldots,T_N\}$.
Then for every
$\epsilon>0$ the Zariski closure of the set of solutions $(u,v)\in\Gamma$ 
of the inequality
$$
h(f(u,v))<(1-\epsilon)\max\{h(T_1(u,v)),\ldots,h(T_N(u,v))\}\eqno(1.4)
$$
is a finite union of translates of proper subtori of $\Gm$.
Also, outside a finite union of translates of proper subtori of $\Gm$,
one has
$$
 h(f(u,v))>(1-\epsilon) 
\max\left\{{h(u)\over 2\deg_Y f },
{h(v)\over 2\deg_X f }\right\}.\eqno(1.5)
$$
}
           
\bigskip

We shall now explain the strategy behind our proofs.
For instance, to prove our Main Theorem above,
observe first that we can reduce to the case when $f$ is defined 
over a number field $K$  and $\Gamma$ is the group 
$(\OS^\times)^2 $ of points  whose coordinates are $S$-units
for a certain fixed finite set of places $S\subset M_K$. Now,
put $f(X,Y)={p(X,Y)\over q(X,Y)}$, where $p,q\in K[X,Y]$ 
are coprime polynomials;  for every point
$(u,v)\in K^2$, where $f$ is defined, the height of $f(u,v)$
is $h(f(u,v))=h(p(u,v):q(u,v))$. Hence it verifies
$$
\eqalign{
h(f(u,v))&=\sum_{\mu\in M_K}\log\max\{|p(u,v)|_\mu,|q(u,v)|_\mu\}\cr
&=\sum_{\mu\in M_K}\log^+ \max\{|p(u,v)|_\mu,|q(u,v)|_\mu\}+
\sum_{\mu\in M_K}\log^- \max\{|p(u,v)|_\mu,|q(u,v)|_\mu\}
}
$$
where $\log^+(\cdot)=\max\{0,\log(\cdot)\}$ and $\log^-(\cdot)=
\min\{0,\log(\cdot)\}$. Hence
$$
h(f(u,v))\geq\max\{h(p(u,v)),h(q(u,v))\}-
\sum_{\mu\in M_K}
\log^- \max\{|p(u,v)|_\mu,|q(u,v)|_\mu\}.
$$
We then need   lower bounds for  $\max\{h(p(u,v)),h(q(u,v))\}$
and  for   $\sum_{\mu\in M_K}
\log^- \max\{|p(u,v)|_\mu,|q(u,v)|_\mu\}$. These estimates are the object
of Propositions 1 and 4 and will be proved by using in an essential way the
 Subspace Theorem. While the techniques to prove lower bounds like   in
Proposition 1 are well known, to prove upper bounds as in Proposition 4
we need a new method, introduced in [5] and developed in
[1] and [6], which consists in applying the Subspace Theorem  
to suitable linear combinations of $S$-units and $S$-integers.\medskip

\noindent{\bf Acknowledgements.} After the present paper was written, Prof. J. Silverman
kindly sent us a preprint exploring related questions. Concerning the formulations in terms
of heights, he goes further than us, relating  with  certain cases of  Vojta's conjecture:
for instance he proves that Corollary 1 follows from the Vojta conjecture suitably applied to
the  blow-up of $\Pr_1^2$ at one point. A paper of Silverman in this respect is expected soon.

The authors are grateful to the referee for his comments.

\bigskip

\noindent{\bf \S 2 Proofs}. We start by recalling the definitions:
$$
\log^+x:=\max \{0,\log x\},\qquad \log^-x:=\min\{0,\log x\}\qquad x>0.
$$
Note that these functions are nondecreasing.

In the sequel, $K$ will denote a number field, $M_K$
the set of places of $K$, $M_0$ the subset of finite places;
for a finite subset $S\subset M_K$ containing
all archimedean places, let us denote by $\OS$ the ring of $S$-integers
and by $\OS^\times$ its group of units ($S$-units).
For a place $\mu\in M_K$, let us denote by $|\cdot|_\mu$ the corresponding
absolute value, normalized {\it with respect to $K$}, i.e. in such a way
that the product formula holds, and the absolute Weil height reads 
$H(x)=\exp(h(x))=\prod_\mu \max\{1,|x|_\mu\}$. 
For a point $\x=(x_1,\ldots ,x_n)\in K^n$
and a place $\mu\in M_K$, put $|\x|_\mu=\max_i |x_i|_\mu$.
For a point 
$\x=(x_1:\ldots :x_n)\in \Pr_{n-1}(K)$, 
put $H(\x)=\prod_\mu \max\{|x_1|_\mu,\ldots,|x_n|_\mu\}=\prod_\mu |\x|_\mu$ 
(where in the last term the point $\x$ is considered in $K^n$), and call
it the height of the point $\x$; note that it is well defined
in view of the product formula.
\bigskip

We beging by stating the  so called Subspace Theorem, 
which represents the main technical tool in the present proofs.
\medskip

\noindent {\bf Subspace Theorem}. 
{\it Let $S$ be a finite set of absolute values
of a number field $K$, 
including the archimedean ones. For $\mu\in S$, let $L_{1,\mu},
\ldots,L_{N,\mu}$ be linearly independent linear forms in $N$ variables,
defined over $K$; let $\delta>0$. Then the solution 
$P=(x_1,\ldots,x_N)\in K^N\setminus\{0\}$ to the inequality 
$$
\prod_{\mu\in S}\prod_{i=1}^N{|L_{j,\mu}(P)|_\mu\over |P|_\mu}<
H(P)^{-N-\delta}
$$
are contained in finitely many proper subspaces of $K^N$.}
\medskip

For a proof, {\it see} [11,12].

\medskip 

In the sequel, we shall also make use of the following two lemmas:
\medskip

\noindent{\bf Lemma 1}. {\it Let 
$\tilde{K}$ be a number field, $\varphi(T)\in \tilde{K}(T)$ 
be a rational function of degree $d$. Then for $t\in \tilde{K}$,
$$
H(\varphi(t))\gg H(t)^d,
$$
where the constant appearing implicitely in the symbol $\gg$ can
be effectively determined.}
\smallskip

The proof of Lemma 1 can be found for instance in [9]: it is a very
particular case of Theorem B.2.5 of [9].
\medskip

\noindent{\bf Lemma 2}. {\it  Let $\{T_1(X,Y),\ldots,T_M(X,Y)\}$ be a
set of non constant monomials in $X,Y$ and suppose that $\deg_X T_i\leq d_1,
\deg_Y T_i \leq d_2$ for $i=1,\dots,M$. Suppose also that the intersection
of the $M$ subgroups defined by the equations
$T_i(X,Y)=1$ is zero-dimensional.
Let $(u,v)\in (K^\times)^2$ be an algebraic point  of $\Gm$. Then 
$$
\max\{h(T_1(u,v)),\ldots,h(T_M(u,v))\}\geq {1\over 2} \max\left\{{h(u)\over d_2},
{h(v)\over d_1}\right\}.
$$
}

\noindent{\it Proof}. 
Denote by $T\subset\Z^2$  the set of exponents
$(i,j)$ appearing in the monomials $T_1,\ldots,T_M$.
The group generated by $T$ in $\Z^2$ is of rank two by our assumption.
Hence the set $T$ contains two linearly independent  points
$(i_1,i_2), (j_1,j_2)$.
Hence  $X^d=(X^{i_1}Y^{i_2})^{j_2}(X^{j_1}Y^{j_2})^{-i_2}$
where $d=i_1j_2-i_2j_1\neq 0$.
By substituting $X,Y$ by $u,v$ and considering the corresponding 
height one obtains
$$
|d|\cdot h(u)\leq |j_2|\cdot h(u^{i_1}v^{i_2})+|i_2|\cdot h(u^{j_1}v^{j_2})
\leq 2 d_2\max\{h(u^{i_1}v^{i_2}),h(u^{j_1}v^{j_2})\}.
$$
Hence the maximum of the heights of the monomials 
is bounded from below by 
${1\over 2} {h(u)\over d_2}$. Of course the analogous 
inequality holds for $h(v)$, thus proving the lemma.
\medskip

The following proposition is essentially a corollary of  Evertse's paper
[8]. For the readers'
convenience, we give here a complete proof.
\medskip

\noindent{\bf Proposition 1}. {\it Let $S\subset M_K$ and 
$\Gamma\subset \Gm(\barQ)$ be as above,
$f(X,Y)\in K [X,Y,X^{-1},Y^{-1}]$ 
be a regular function on $\Gm$. Denote by $T_1(X,Y),\ldots,T_N(X,Y)$
the monomials in $X,X^{-1},Y,Y^{-1}$ which appear in $f$
and suppose that $1\in\{T_1,\ldots,T_N\}$.
Let $\epsilon>0$ be fixed. Then the solutions $(u,v)\in\Gamma$
to the inequality
$$
 h(f(u,v))<(1-\epsilon) \max\{h(T_1(u,v)),\ldots,h(T_N(u,v))\}\eqno(2.1)
$$
lie in the union of finitely many translates of
1-dimensional  subtori of $\Gm$, which can be effectively determined, and a
finite set. The same conclusion holds for the solutions to the inequality
$$
\sum_{\mu\in S}\log^-|f(u,v)|_\mu<-\epsilon\max\{h(u),h(v)\}.\eqno(2.2)
$$
}
\medskip

The idea of the proof is two-fold: first, by an application of the 
Subspace Theorem, we obtain that the solutions to the above inequality
belong to the union of finitely many proper subgroups. Then, a Liouville-type
argument enables to bound the degree of the minimal equations satisfied by  any solution.
 
\medskip

\noindent{\it Proof of Proposition 1}. We shall first prove our contention
for  the solutions to the inequality (2.1). 
Let us write 
$$
f(X,Y)=\sum_{(i,j)\in T}a_{i,j}X^iY^j=b_1T_1+\ldots+b_NT_N
$$
where $T\subset\Z^2$ is a finite set of cardinality $N$ and,
for $(i,j)\in T$, $a_{i,j}\in K^\times$  
are nonzero algebraic numbers in $K$ ordered as
$b_1,\ldots,b_N$. 

We first treat the (easier) case where $T$ is contained in a line of
$\Z^2$. In this case, there exists a rational function $\varphi(T)\in K[T,{1\over T}]$
and a pair $(i_0,j_0)$ such that identically
$$
f(u,v)=\varphi(u^{i_0}v^{j_0}).
$$
Also, for $(u,v)\in\Gm(K)$, $\deg(\varphi)\cdot h(u^{i_0}v^{j_0})\geq
\max\{h(T_1(u,v)),\ldots,h(T_N(u,v))\}$. So, by Lemma 1, there exists a computable
constant $c$ such that for all $(u,v)\in\Gm(K)$ 
$$
h(f(u,v))\geq \max\{h(T_1(u,v)),\ldots,h(T_N(u,v))\}-c.
$$
Then, if inequality (2.1) holds, the above lower bound  gives
$$
\max\{h(T_1(u,v)),\ldots,h(T_N(u,v))\}\leq c/\epsilon.
$$
This implies that the pair $(u,v)$ is contained in a finite translate of
a subtorus of $\Gm$ defined by an equation of the form
$$
u^{i_0}v^{j_0}=\lambda
$$
where $\lambda\in K$ is an algebraic number with $h(\lambda)\leq c/\epsilon$.
This implies our contention in that particular case.
 
We now consider  the generic case where the subset $T\subset \Z^2$ is not
contained in a line. 

As we have already remarked, we can choose a finite
set $S$ of places of $K$ containing the archimedean ones 
such that the subgroup $\Gamma$ is contained in $(\OS^\times)^2$. 
Now, for each solution $(u,v)\in (\OS^\times)^2$ of the inequality
(2.1) and each $\mu\in S$, let $k(u,v;\mu)$ be an index
in $\{1,\ldots,N\}$ for which the quantity 
$|T_k(u,v)|_\mu$ is maximal. 

By partitioning the set  of solutions of inequality (2.1)
in at most $N$ subsets, and working separately with each subset, 
we can suppose that for each $\mu\in S$,
$k(u,v;\mu)=:k_\mu$ does not depend on $(u,v)$. 
 
Let us define, for each $\mu\in S$, $N$ linearly independent linear forms
in $T_1,\ldots,T_N$ as follows:   for $i\neq k_\mu$, put 
$L_{i,\mu}=T_i$, while for $i=k_\mu$, put 
$$
L_{k_\mu,\mu}=b_1T_1+\ldots+b_N T_N.
$$
Note that for each $\mu\in S$, the above defined linear forms are indeed
independent, since all coefficients $b_i$ are nonzero.

For a solution $(u,v)\in(\OS^\times)^2$ of the inequality (2.1) we put
$P=P(u,v)=(T_1(u,v),\ldots,T_N(u,v))\in\OS^N$, so that for each $\mu\in S$ we have
$L_{k_\mu,\mu}(P)=f(u,v)$. Note that, by Lemma 2 and the assumption that
$T$ is not contained in a line, the height of $P$ tends to infinity on
every infinite sequence of points $(u,v)\in K^2$.
We now consider the double product 
$$
\prod_{\mu\in S}\prod_{i=1}^N {|L_{i,\mu}(P)|_\mu\over |P|_\mu};\eqno(2.3)
$$
an application of the Subspace Theorem to the above double product will give
what we want. We note that, since the coordinates of $P$ are $S$-integers
(in fact even $S$-units), the quantity $\prod_{\mu\in S}|P|_\mu$ equals
the height $H(P)$ of $P$. We can then rewrite
the above double product as
$$
\prod_{\mu\in S}\prod_{i=1}^N {|L_{i,\mu}(P)|_\mu\over |P|_\mu}
=H(P)^{-N}\cdot (\prod_{\mu\in S} |f(u,v)|_\mu)\cdot
(\prod_{\mu\in S}\prod_{i=1}^N |T_i(u,v)|_\mu)\cdot
\prod_{\mu\in S}|T_{k_\mu,\mu}(u,v)^{-1}|_\mu.\eqno(2.4)
$$ 
Clearly, the inequality $\prod_{\mu\in S} |f(u,v)|_\mu\leq H(f(u,v))$ holds and
by assumption (2.1) we also have $H(f(u,v))<\max\{H(T_1(u,v)),\ldots,
H(T_N(u,v))\}^{1-\epsilon}$. In turn, since by assumption
$1\in\{T_1,\ldots,T_N\}$, the height of $P$ is larger then the 
maximum of the heights of its coordinates. Then from (2.1) it follows that
$H(f(u,v))<H(P)^{1-\epsilon}$. 

Also the third factor in (2.4) equals $1$ since all the $T_i(u,v)$ are $S$-units.
The last factor $\prod_{\mu\in S}|T_{k_\mu,\mu}(u,v)|_\mu$ is equal to the 
height $H(P)$, by the choice of the $k_\mu$. 
Hence the double product (2.3) can be bounded as
$$
\prod_{\mu\in S}\prod_{i=1}^N {|L_{i,\mu}(P)|_\mu\over |P|_\mu}\leq
H(P)^{-N-\epsilon}.
$$
The Subspace Theorem, applied with $\delta=\epsilon$, asserts that
the points $P$ 
arising from the solutions $(u,v)$ to the inequality (2.1) 
lie in finitely many hyperplanes. Then we obtain finitely many
equations  of the form
$\psi_j(u,v)=0$ ($j=1,\ldots,M$), such that:
(i)
$\psi_j(X,Y)$ ($j=1,\ldots,M$)
are irreducible polynomials in $\bar{K}[X,Y]$;
(ii) all but finitely many
solutions of (2.1) satisfy one of these equations and (iii) each such equation
has infinitely many solutions in $S$-units $(u,v)\in(\OS^\times)^2$.
By a Theorem of Lang [10, Tm. 7.3, p. 207], 
each equation $\psi_j(X,Y)=0$ defines a 
translate of a proper subtorus.
Every translate of a proper subtorus
is either a point or a curve defined by an equation of the form
$$
U^pV^q=w\eqno(2.5)
$$
for coprime $p,q$ and nonzero $w\in K^\times$. 
Our next goal is to show that
the pairs $(p,q)$ such that infinitely many solutions of $(2.1)$ satisfy the
above equation (2.5) can be effectively determined.

Fix such a translate, given by (2.5), and suppose it contains
infinitely many solutions $(u,v)\in(\OS^\times)^2$ to the inequality (2.1).
For each $(u,v)\in (\OS^\times)^2$, write
$$
u=t^q,\qquad v=s^{-p}
$$
where $t,s$ lie in a fixed number field $\tilde{K}$. 
The relation $u^pv^q=w$
can be rewritten as $(t/s)^{pq}=w$; so that $(t/s)^p$ is  
a $q$-th root of $w$. Select one such root $\overline{w}$. 
We obtain the parametrization
$$
\left\{\eqalign{
u&=t^q\cr
v&=\overline{w}t^{-p}
}\right.\eqno(2.6)
$$

Then, using the above parametrization,  we can write 
$$
f(u,v)=f(t^q,\overline{w}t^{-p})=:\varphi(t),\eqno(2.7)
$$
for a rational function $\varphi(t)\in \tilde{K}(t)$.
 
Note that in our case  the function $\varphi(t)$ in (2.7) can be written as
$$
\varphi(t)=\sum_{(i,j)\in T}a_{i,j}\overline{w}^jt^{qi-pj}=\sum_{l\in\Z} c_lt^l
$$
where, for every $l\in\Z$, 
$$
c_l=\sum_{(i,j)\in T\ :\ qi-pj=l}a_{i,j}\overline{w}^j,
$$
the last sum running over those pairs $(i,j)\in T$ with $qi-pj=l$.
Since the degree $d$ of $\varphi(t)$ is at least
$\max\{|l| : c_l\neq 0\}$,  Lemma 1 implies  that 
$$
h(f(u,v))=h(\varphi(t))\geq \max\{|l| : c_l\neq 0\}h(t) +O(1)
$$ 
Now, if $l=qi-pj$, for a pair $(i,j)\in T$, the height of
the corresponding monomial $u^iv^j$ is $|l|h(t)+O(1)$. Hence,
if all the coefficients $c_l$, for $l$ of the form
$l=qi-pj$ for at least one  $(i,j)\in T$, are
nonzero, we have the inequality
$$
h(f(u,v))\geq \max\{h(T_1(u,v)),\ldots,h(T_N(u,v))\}+O(1)
$$
which contradicts (2.1) for large values of $\max\{h(T_1(u,v)),\ldots,h(T_N(u,v))\}$
(recall that by Lemma 2 such a maximum does tend to infinity, so 
we indeed obtain a contradiction with the hypothesis that the given
subgroups contains infinitely many solutions of (2.1)). 
Hence some $c_l$ must vanish, which means that,
after the substitution (2.6), some cancellation occurs among the monomials 
$T_1,\ldots,T_N$.  This fact implies that there exist two distinct pairs
of exponents $(i,j),(i^\prime,j^\prime)\in T$ such that
$qi-pj=qi^\prime-pj^\prime$, so $(p,q)$ verifies a (non trivial) equation of the form
$$
q(i-i^\prime)-p(j-j^\prime)=0.
$$
Each such equation has exactely two solutions in coprime integers $p,q\in\Z$, which
can be effectively determined. This proves the first contention in Proposition 1.

We shall now obtain the same conclusion under assumption (2.2). Chosing linear forms
as before, we arrive again at equality (2.4). 
As we have already remarked, the third factor
equals $1$, since $T_i(u,v)$ are $S$-units for all $i=1,\ldots,N$.
The product of the second and fourth factors in (2.4) can be rewritten
as
$$
\eqalign{(\prod_{\mu\in S} |f(u,v)|_\mu) 
\prod_{\mu\in S}|T_{k_\mu,\mu}(u,v)^{-1}|_\mu & =
\left(\prod_{\mu\in S} \min\{1,|f(u,v)|_\mu\}\right)\cdot
\cr &\left(\prod_{\mu\in S} \max\{1,|f(u,v|_\mu)\}\right)
\cdot \prod_{\mu\in S}|T_{k_\mu,\mu}(u,v)^{-1}|_\mu.}
$$
Clearly, for every $\mu\in S$, 
$$
\max\{1,|f(u,v)|_\mu\}\leq \max\{1,|N|_\mu\}\cdot
\max\{|b_1|_\mu,\ldots,||b_N|_\mu\}\cdot |T_{k_\mu,\mu}(u,v)|_\mu,
$$
so  we have
$$
\left(\prod_{\mu\in S} \max\{1,|f(u,v)|_\mu\}\right) 
\prod_{\mu\in S}|T_{k_\mu,\mu}(u,v)^{-1}|_\mu=O(1).
$$
Now from (2.4) we obtain (recall that the third factor in (2.4) equals $1$)
$$
\prod_{\mu\in S}\prod_{i=1}^N {|L_{i,\mu}(P)|_\mu\over |P|_\mu}
\ll H(P)^{-N}\cdot \prod_{\mu\in S} \min\{1,|f(u,v)|_\mu)\}
$$
where the constant implicit in the symbol $\ll$ does not depend on $(u,v)$.
The inequality (2.2) reads 
$$
\prod_{\mu\in S}
\min\{1,|f(u,v)|_\mu)\}<\max\{H(u),H(v)\}^{-\epsilon}.
$$
Clearly, 
$H(P)\leq \max\{H(u),H(v)\}^{D}$ for a positive constant $D$
(which can be taken to be $2\max\{|i|+|j|\, :\, (i,j)\in T\}$).
Then from the last displayed inequality we obtain
$$
\prod_{\mu\in S}\prod_{i=1}^N {|L_{i,\mu}(P)|_\mu\over |P|_\mu}
\ll H(P)^{-N-\epsilon/D}.
$$
The Subspace Theorem, applied with $\delta=\epsilon/(2D)$, say, 
asserts that all points $P(u,v)$ are contained in the union of
finitely many hyperplanes. The final argument to obtain the conclusion
is exactely the same as in the previous case, so we do not repeat it.

\bigskip

Our next goal is to prove an ``explicit"  version of Corollary 1,
which will be used in the subsequent proofs of Theorem 1 and the Main Theorem.
Corollary 1 can be refrased by saying that if the ratio
$(u-1)/(v-1)$ has ``small" height, where $u$ and $v$ belong to a given
finitely generated multiplicative group, then, apart finitely many exceptions,
the pair $(u,v)$ verifies a multiplicative dependence relation, i.e.
a relation of the form 
$$
u^p\cdot v^q=1
$$
for a vector $(p,q)\in\Z^2\setminus\{(0,0)\}$. In other words,
the point $(u,v)\in\Gm$ belongs to a one-dimensional subgroup.
Note that if $u$ and $v$ are not both roots of unity,
then  vectors $(p,q)\in\Z^2$ such that the above relation holds
form a rank-one subgroup. 
The following Proposition 2 quantifies
the above mentioned Corollary, in the sense that it permits to explicitely
obtain the dependence relations satisfied by all but finitely many
exceptions to the inequality (1.3).
 
\medskip 

\noindent{\bf Proposition 2}. {\it  Let $K$ be a number field, 
$S$ a finite set of places as before and $\epsilon>0$ a real
number.  All but finitely many  solutions $(u,v)\in(\OS^\times)^2$
to the inequality
$$
\sum_{\mu\in M_K} \log^-\max\{|u-1|_\mu, |v-1|_\mu\}<-\epsilon
\cdot \max\{h(u),h(v)\}
\eqno(2.8)
$$
are contained in finitely many 1-dimensional subgroups of
$\Gm$ which can be effectively determined. Namely they are defined by
an equation of the form $u^p=v^q$ for coprime integers $p,q$ with
$\max\{|p|,|q|\}\leq \epsilon^{-1}$.
}

\medskip

\noindent {\it Proof of Proposition 2}. 
Let us suppose that $(u_i,v_i)$ is an infinite sequence
in $(\OS^\times)^2$ verifying the above inequality (2.8). 
We shall prove the existence of finitely many pairs
$(p,q)\in\Z^2\setminus\{(0,0)    \}$ such that for all but finitely many indices,
$u_i$ and $v_i$ verifie one dependence relation of the form
$u^pv^q=1$ (in particular for all  but finitely many indices, $u_i,v_i$
are  multiplicatively dependent). We shall later show how to bound
the size of $p$ and $q$. In the sequel, we will drop the index $i$ 
for convenience of notation.

First of all, notice that by Ridout's theorem (the case $N=2$ of the 
above Subspace Theorem) for all but finitely many $(u,v)\in\Gamma$
we have
$$
\prod_{\mu\in S}\min\{1,|u-1|_\mu\}>H(u)^{-\epsilon},\qquad
{\rm and}\qquad 
\prod_{\mu\in S}\min\{1,|v-1|_\mu\}>H(v)^{-\epsilon}.
$$
Hence we clearly have that apart finitely many pairs $(u,v)\in(\OS^\times)^2$
$$
\sum_{\mu\in S} \log^-\max\{|u-1|_\mu, |v-1|_\mu\}>-{\epsilon\over 2}
\cdot \max\{h(u),h(v)\};
$$
so  all but finitely many solutions of (2.8)  
also verify 
$$
\sum_{\mu\in M_K\setminus S} \log^-\max\{|u-1|_\mu, |v-1|_\mu\}
=\sum_{\mu\in M_K\setminus S} \log \max\{|u-1|_\mu, |v-1|_\mu\}<
-{\epsilon\over 2} \cdot \max\{h(u),h(v)\}. \eqno(2.8')
$$
We can suppose (since the expressions (2.8) , (2.8') are symmetric in $u,v$) 
that $H(v)\geq H(u)$. 
Let us denote by $S^+$ (resp. $S^-$), the subset of $S$ made up of those
 absolute values $\mu$ such that $|v|_\mu>1$ (resp. $|v|_\mu\leq 1$). 
Since there are only finitely many choices  for the pair
$S^+,S^-$, we can partition the solutions in a finite number of subsets for which
$S^+$ and $S^-$ are constant. From now on, we shall work separately with each
subset.

Now the proof follows the same lines as [6].
Let $(u,v)$ be a pair  in $(\OS^\times)^2$, satisfying inequality (2.8').
We write, for a positive integer $j$,
$$
z_j(u,v)=z_j:={u^j-1\over v-1}.\eqno(2.9)
$$
Note that 
$$
z_j=z_1\cdot(u^{j-1}+\ldots+u+1).\eqno(2.10)
$$
Fix an integer $h>0$  and consider the identity
$$
{1\over v-1}={1\over v}\cdot {1\over 1-v^{-1}}={1\over v}\left(1+v^{-1}+\ldots
+v^{-h+1}+{v^{-h}\over 1-v^{-1}}\right).
$$
Let us fix a second integer $k>0$.
Then for  $j\in\{1,\ldots,k\}$ 
we obtain, on multiplying by $u^{j}-1$ in the above identity,
$$
z_j=(u^j-1)\cdot
\left(v^{-1}+v^{-2}+\ldots+v^{-h}+{v^{-h-1}\over v-1}\right).
$$
For $\mu\in S^+$, we then derive the inequality (in fact an equality)
$$
|z_j-u^jv^{-1}-\ldots-u^jv^{-h}+v^{-1}+\dots+v^{-h}|_\mu \leq
|u^j-1|_\mu \cdot {|v|_\mu^{-h-1}\over |1-v|_\mu}.\eqno(2.11)
$$
We put $N=hk+h+k$; for convenience we shall write vectors in $K^N$ as
$$
\x=(x_1,\ldots,x_N)=(z_1,\ldots,z_k,y_{0,1},\ldots,y_{0,h},\ldots,
y_{k,1},\ldots,y_{k,h}).
$$
In this notation we choose linear forms with integral coefficients
as follows. For $j=1,\ldots,k$, and $\mu\in S^+$, we put
$$
L_{j,\mu}=z_j+y_{0,1}+\ldots y_{0,h}-y_{j,1}-\ldots-y_{j,h},
$$
 while, for $(j,\mu)\not\in\{1,\ldots,k\}\times S^+$ we put
$$
L_{j,\mu}=x_j.
$$
Observe that for each $\mu\in S$, the linear forms are indeed linearly
independent.
For a given pair $(u,v)$ we also set
$$
P=P(u,v)=(z_1,\ldots,z_k,v^{-1},\ldots,v^{-h},uv^{-1},\ldots,
uv^{-h},\ldots,u^kv^{-1}\ldots,u^kv^{-h}).
$$
Note that all the coordinates of $P$ apart the first $k$ of them
are $S$-units.

To apply the Subspace Theorem we have to estimate the  
double product 
$\prod_{j=1}^N\prod_{\mu\in S}{|L_{j,\mu}(P)|_\mu\over |P|_\mu}$
at points $P=P(u,v)$.
First note that for $j>k$, 
$$
\prod_{\mu\in S}|L_{j,\mu}(P)|_\mu=1
$$
since all the involved linear forms are projections on the coordinates,
which are $S$-units.
Then
$$
\prod_{j=k+1}^N\prod_{\mu\in S}{|L_{j,\mu}(P)|_\mu\over |P|_\mu}=
\prod_{\mu\in S}|P|_\mu^{-N+k}.
\eqno(2.12)
$$
Let now $j$ be an index in $\{1,\ldots,k\}$. For $\mu\in   S^-$, we
recall that $L_{j,\mu}(P)=z_j=(u^j-1)/(v-1)$, so in particular
$$
{|L_{j,\mu}(P)|_\mu\over |P|_\mu}\leq {|u^j-1|_\mu\over |v-1|_\mu}
\cdot |P|_\mu^{-1}.\eqno(2.13) 
$$
For $\mu\in S^+$ we use $(2.11)$ which gives, since $ L_{j,\mu}(P)$ coincides
with the left side term of $(2.11)$,
$$
|L_{j,\mu}(P)|_\mu
\leq  |u^j-1|_\mu \cdot {|v|_\mu^{-h}\over |1-v|_\mu}
$$
so, since $|u^j-1|_\mu\leq \max\{1,|2|_\mu\}\cdot \max\{1,|u|_\mu\}^j$,
$$
{|L_{j,\mu}(P)|_\mu\over |P|_\mu}\leq 
\max\{1,|2|_\mu\}\cdot \max\{1,|u|_\mu\}^j\cdot {1\over |1-v|_\mu}\cdot
 |v|_\mu^{-h}\cdot |P|_\mu^{-1}. \eqno(2.14)
$$
Then, from (2.13), (2.14), we have
$$
\prod_{j=1}^k \prod_{\mu\in S}{|L_{j,\mu}(P)|_\mu\over |P|_\mu}\leq
\left(\prod_{j=1}^k \prod_{\mu\in S}
\max\{1,|2|_\mu\}\cdot \max\{1,|u|_\mu\}^j\cdot {1\over |1-v|_\mu}
\cdot |P|_\mu^{-1}
\right)\cdot\left(\prod_{\mu\in S^+}|v|_\mu^{-h}\right)^k.
$$
Now, using the fact that $\prod_{\mu\in S} \max\{1,|2|_\mu\}\leq H(2)=2$ 
(and the analogue estimates for the products of $\max\{1,|u|_\mu\}$ and 
${1\over |1-v|_\mu}$), and using also the equality
 $H(v)=\prod_{\mu\in S^+} |v|_\mu$, we obtain
$$
\prod_{j=1}^k\prod_{\mu\in S}
{|L_{j,\mu}(P)|_\mu\over |P|_\mu}\leq 
H(v)^{-hk}(2H(u))^{1+\ldots+k} H(1-v)^k\prod_{\mu\in S}|P|_\mu^{-k}.
\eqno(2.15)
$$
Finally, from (2.12), (2.15) we obtain 
$$
\prod_{j=1}^N\prod_{\mu\in S} {|L_{j,\mu}(P)|_\mu\over |P|_\mu}
\leq \left(\prod_{\mu\in S}|P|_\mu^{-N+k}\right)  H(v)^{-hk}
(2H(u))^{k^2}
H(1-v)^k \prod_{\mu\in S}|P|_\mu^{-k}.
$$       
Observe that the product of the first and last factors can be written as
$$
\left(\prod_{\mu\in S}|P|_\mu\right)^{-N}=H(P)^{-N}\cdot
\left(\prod_{\mu\not\in S}|P|_\mu\right)^N.
$$
Hence we can rewrite the above inequality as
$$
\prod_{j=1}^N\prod_{\mu\in S} {|L_{j,\mu}(P)|_\mu\over |P|_\mu} \leq
H(P)^{-N} \left(\prod_{\mu\not\in S}|P|_\mu\right)^N
 H(v)^{-hk} (2H(u))^{k^2} H(1-v)^k \eqno(2.16)
$$
Our next goal will be to estimate the quantity $\prod_{\mu\not\in S}|P|_\mu$.
 Now, observe that the only
coordinates which are not $S$-integers are the first $k$, i.e.
$z_1,\ldots,z_k$, and, for $j=1,\ldots,k$,
$$
z_j={u^j-1\over v-1}=z_1\cdot (u^{j-1}+u^{j-2}+\ldots +u+1);
$$
also, the factor $u^{j-1}+\ldots+u+1$ is an $S$-integer, 
so the only contribution comes from the first factor $z_1=(u-1)/(v-1)$;
hence  we have the bound
$$
\prod_{\mu\not\in S}|P|_\mu\leq   
\prod_{\mu\not\in S}\max\left\{1,\left|{u-1 \over v-1}\right|_\mu\right\}. 
$$
Now, we have
$$
\prod_{\mu\not\in S}\max\left\{1,\left|{u-1 \over v-1}\right|_\mu\right\}=
\prod_{\mu\not\in S}{1\over |v-1|_\mu}\cdot \max\{|u-1|_\mu, |v-1|_\mu\}
\leq H(v-1)\prod_{\mu\not\in S}\max\{|u-1|_\mu, |v-1|_\mu\};
$$
we now use the hypothesis that our pairs $(u,v)$ satisfy (2.8')
so that 
$$
\prod_{\mu\not\in S}\max\{|u-1|_\mu, |v-1|_\mu\}<H(v)^{-\epsilon/2},
$$
so finally 
$$
\prod_{\mu\not\in S}|P|_\mu\leq   H(v-1)\cdot H(v)^{-\epsilon/2}.
$$
Using the above estimate and the fact that $H(1-v)\leq 2H(v)$,  
 we can rewrite (2.16) as 
$$
\eqalign{
\prod_{j=1}^N\prod_{\mu\in S} {|L_{j,\mu}(P)|_\mu\over |P|_\mu} &\leq
H(v)^{-hk}
H(P)^{-N} H(v)^{(1-\epsilon/2)N} 2^N (2H(u))^{k^2} (2H(v))^{k}\cr
&\leq H(v)^{-hk}
H(P)^{-N} H(v)^{(1-\epsilon/2)N} 2^N (2H(v))^{k^2+k}. 
} \eqno(2.17)
$$
Since $N=hk+h+k$, the exponent of $H(v)$ can be estimate (at least
when $k\geq 2$), by
$$
\left(1-{\epsilon\over 2}\right)N+k^2+k-hk<-{\epsilon\over 2} hk + h+k+k^2+k\leq
-{\epsilon\over 2} hk+  h+2k^2.
$$
By choosing $k>{4\over \epsilon}$ we obtain that 
 for large enough $h$ (in particular $h>2k^2+1$), 
 ${\epsilon hk/2}-h-2k^2=:\epsilon_0>0$. 
Choose $k,h$ in such a way. 
From (2.17) we then obtain that the double product is bounded as
$$
\prod_{j=1}^N\prod_{\mu\in S} {|L_{j,\mu}(P)|_\mu\over |P|_\mu} \leq
H(P)^{-N}  H(v)^{-\epsilon_0}\cdot 2^{N+k+k^2}.  \eqno(2.18)
$$
Our next goal, in view of the application of the Subspace Theorem,
will be to compare the height of $P$ with the height of $v$.
A rough estimate gives
$$
H(P)\leq H(u)^k H(v)^h H(1-v)\leq H(u)^k 2 H(v)^{h +1}\leq
H(v)^{h+k+2}
$$
at least for $H(v)\geq 2$, which we may suppose. Hence
$$
H(v)^{-\epsilon_0}\leq H(P)^{{-\epsilon_0\over h+k+2}}
$$
so that, for $\delta:={\epsilon_0\over h+k+3}$ and large values of
$H(v)$, we obtain from (2.18) 
$$
\prod_{j=1}^N\prod_{\mu\in S} {|L_{j,\mu}(P)|_\mu\over |P|_\mu}
\leq H(P)^{-N-\delta}
$$    
Now we are able to apply the Subspace Theorem, which implies
there exist  finitely many  hyperplanes containing all the points
$P$. Choose one such  hyperplane, containing infinitely
many points $P$. 
Then we obtain an equation of the kind
$$
\eta_k{u^k-1\over v-1}+\ldots+\eta_{1} {u-1\over v-1}+
\sum_{i,j} \rho_{i,j} u^jv^i=0, \eqno(2.19)
$$
where the sum runs over the pairs $(i,j)\in\{-1,\ldots,-h\}\times\{0,\ldots,k\}$. 
Here $\eta_1,\ldots,\eta_k$ and the $\rho_{i,j}$ are
elements of the number field $K$, not all zero.
By multiplying by $v^h(v-1)$ we obtain that the point $(u,v)$ lies
in the affine curve given by the equation
$$
\eta_k V^h(U^k-1)+\ldots +\eta_1 V^h(U-1)+
(V-1)\sum_{i,j} \rho_{i,j}U^jV^{h+i}=0. 
$$
We pause to show that this equation is non trivial, 
so it defines in fact a curve.
If the left side term vanishes identically, then the
binomial $(V-1)$ would divide $\eta_k V^h(U^k-1)+\ldots +\eta_1 V^h(U-1)$,
which implies that $\eta_1=\ldots=\eta_k=0$. But then all coefficients
$\rho_{i,j}$ would vanish, contrary to the assumption.
    
In view of (2.19), infinitely many points $(u,v)$
lie in the curve defined by the above equation. Then, 
by the mentioned  theorem of Lang  [L, Thm. 7.3, p.207],
all such points lie in the union of finitely many translates
of subtori of ${\bf G}_m^2$;  every translate of a proper subtorus
is either a point or a curve defined by an equation of the form (2.5) 
for coprime $p,q$. The remaining of the proof is similar to the last part
of the proof of Proposition 1:
Fix such a translate and suppose it contains
infinitely many solutions $(u,v)\in(\OS^\times)^2$ to the inequality (2.8).
For each $(u,v)\in (\OS^\times)^2$, write
$$
u=t^q,\qquad v=s^{-p}
$$
where $t,s$ lie in a fixed number field $\tilde{K}$. 
As in the proof of Proposition 1, we arrive at the parametrization (2.6).
Then  
$$
{u-1\over v-1}={t^q-1\over \overline{w}t^{-p}-1}=:\varphi(t),\eqno(2.20)
$$
say. 

Now, we distinguish four cases:
\medskip

\noindent{\it First Case}: $pq>0$ and $\overline{w}^q\neq 1$.

In this case the rational function $\varphi$ in (2.20), which
can be written as $\varphi(t)=-t^p\cdot{t^q-1\over t^p-\overline{w}}$,
has degree $|p|+|q|$ (note that the numerator and denominator
are coprime polynomials in $t$, or in $t^{-1}$). 
So, in this case, we have $H((u-1)/(v-1))\gg H(u)\cdot H(v)$, i.e.
$$
h((u-1)/(v-1))\geq h(u)+h(v)+O(1)
$$
where the implied constant in the $O(1)$ term does not depend on $u,v$.
Since
$$
h((u-1)/(v-1))=h(u-1:v-1)\leq h(u)+h(v)+\sum_{\mu\in M_K}\log^-\max\{
|u-1|_\mu,|v-1|_\mu\}+O(1),
$$
we obtain from the above inequalities and from
(2.8) a uniform bound for the $\max\{h(u),h(v)\}$,
contrary to the assumption that the given subgroup contains 
infinitely many solutions to the inequality (2.3).
\medskip

\noindent{\it Second Case:} $pq>0$ and $\overline{w}^q=1$. 

First of all, note that in this case $w=\overline{w}^q=1$, 
so that the dependence relation (2.5) takes the form $u^pv^q=1$.
Also, the rational function $\varphi(t)$ has degree $|p|+|q|-1$, 
since the polynomials $t^p(t^q-1)$, and $t^p-\overline{w}$
have a greatest common divisor of degree one. Then from  
the parametrization (2.6), we clearly obtain
$$
h(1:u-1:v-1)\leq (|p|+|q|)h(t) +O(1)
$$ 
and from Lemma 1
$$
h((u-1)/(v-1)) \geq (|p|+|q|-1)h(t).
$$
Since we always have 
$$
h((u-1)/(v-1))=h(u-1:v-1)\leq h(1:u-1:v-1)
+\sum_{\mu\in M_K}\log^-\max\{|u-1|_\mu,|v-1|_\mu\}
$$
we deduce 
$$
(|p|+|q|-1)h(t)\leq (|p|+|q|)h(t)+\sum_{\mu\in M_K}\log^-\max\{|u-1|_\mu,|v-1|_\mu\}
+O(1)
$$
from which
$$
-\sum_{\mu\in M_K}\log^-\max\{|u-1|_\mu,|v-1|_\mu\}\leq h(t)+O(1).
$$
Since $h(t)=\max\{h(u),h(v)\}(\max\{|p|,|q|\})^{-1}$, the above
inequality can be satisfied by infinitely many solutions of (2.8) only if
$\max\{|p|,|q|\}\leq \epsilon^{-1}$, as wanted.
\medskip

\noindent{\it Third Case}: $pq<0$ and $\overline{w}^q\neq 1$.

Now, the rational function $\varphi$ has degree
$\max\{|p|,|q|\}$, so, by Lemma 1,
we have $H((u-1)/(v-1))\gg \max\{H(u),H(v)\}$, i.e.
$$
h(u-1:v-1)\geq \max\{h(u),h(v)\}+O(1).\eqno(2.21)
$$
Note that in this case it follows from the parametrization
(2.6) that $u$ and $v$ are, apart from a constant $\overline{w}$,
positive powers of a same element $t$ (or $t^{-1}$). From this fact it easily
follows that
$$
h(1:u:v)=\max\{h(u),h(v)\}+O(1). \eqno(2.22) 
$$  
Using again the inequality
$$
h(u-1:v-1)\leq h(1:u:v)+\sum_{\mu\in M_K}\log^-\max\{
|u-1|_\mu,|v-1|_\mu\} 
$$ 
we obtain from (2.21) and (2.22)
$$
\max\{h(u),h(v)\}+O(1)\leq h(u-1:v-1)\leq \max\{h(u),h(v)\} 
+\sum_{\mu\in M_K}\log^-\max\{|u-1|_\mu,|v-1|_\mu\}.
$$
Then from (2.8) we obtain, as before,
a uniform upper bound for the $\max\{h(u),h(v)\}$, contrary to the assumption that
our subgroup contains infinitely many solutions of (2.8).
\medskip

\noindent{\it Fourth Case}:
 $pq<0$ and  $\overline{w}^q=1$.

Then from the parametrization $(2.6)$ we obtain, as in the Second Case,  
$u^pv^q=1$, i.e. a multiplicative dependence relation as wanted. 

Now, since  $p,q$ are coprime and of opposite sign,
the rational function $\varphi(t)=(t^q-1)/(t^{-p}-1)$ has degree
$d=\max\{|q|-1, |p|-1\}$.  Then  the  Lemma 1 gives
$$
h((u-1)/(v-1))\geq \max\{|q|-1,|p|-1\}h(t)+O(1)=(\max\{|q|,|p|\}-1)h(t)+O(1)
$$
On the other hand, the above estimates,  inequality (2.18) (which still
holds), (2.19) and (2.8) give  
$$
h(u-1:v-1)\leq \max\{|q|,|p|\}h(t)-\epsilon \max\{|q|,|p|\}h(t),
$$
so that we finally obtain
$$
(\max\{|q|,|p|\}-1)h(t)+O(1)\leq \max\{|q|,|p|\}h(t)-\epsilon \max\{|q|,|p|\}h(t)
$$
which gives $\max\{|q|,|p|\}\leq \epsilon^{-1}$, as wanted.

Notice that we have also proved that, in the only cases when
infinitely many solutions can occur (i.e. the second and fourth cases), 
the constant $w$ in (2.5) equals $1$.
In other terms, the minimal dependence relation has coprime exponents, 
i.e.  the infinite families of solutions of (2.8) are
contained in finitely many {\it connected } one-dimensional subgroups.
 
\medskip

From Proposition 2 we easily obtain the following statement:
\medskip

\noindent{\bf Proposition 3}. {\it Let $S$, as before,
be a finite set of places of a number field $K$ containing the
archimedean ones,
$\theta,\eta$ be nonzero elements of $K$, $\epsilon>0$. The
solutions $(u,v)\in (\OS^\times)^2$ to the inequality
$$
\sum_{\mu\in M_K}\log^-\max\{|u-\theta|_\mu,|v-\eta|_\mu\}<-
\epsilon\cdot \max\{h(u),h(v)\}
$$
are contained in finitely many translates of proper subgroups of $\Gm$.
Those of dimension one can be effectively determined.
}
\medskip

\noindent {\it Proof}. 
Fix $\theta,\eta\in K^\times$. By enlarging $S$,
we can suppose that they both are $S$-units.
Now, note that  
$$
 \sum_{\mu\in M_K}\log^-\max\{|u-\theta|_\mu,|v-\eta|_\mu\}=
\sum_{\mu\in M_K}\log^-\max\{|u\theta^{-1}-1|_\mu,|v\eta^{-1}-1|_\mu\}+O(1).
$$ 
Then Proposition 2 implies that apart 
finitely many exceptions, the pairs  
 $(u\theta^{-1},v\eta^{-1})$ are contained in finitely many
one-dimensional subgroups, which can be effectivley determined; 
this in turn implies that the pairs
$(u,v)$, apart finitely many of them,
 are contained in finitely many translates of one-dimensional
subgroups, which can be effectively determined.

\medskip

The following Proposition represents a particular but crucial case of Theorem 1.
\medskip 

\noindent {\bf Proposition 4}. {\it Let $r(X),s(X)\in \overline{\Q}[X]$ be 
two non zero polynomials. 
Then for every $\epsilon>0$, all but finitely many
 solutions $(u,v)\in (\OS^\times)^2$ to the inequality
$$
\sum_{\mu\in M_K\setminus S} \log^-\max\{|r(u)|_\mu, |s(v)|_\mu\}
<-\epsilon\cdot (\max\{h(u),h(v)\}) \eqno(2.23)
$$
are contained in finitely many translates of 
one-dimensional subgroups of $\Gm$,  which  
can be effectively determined. Moreover, if   $r(X)$ and $s(X)$ do not both
vanish at $0$, the same conclusion holds for the
solutions to the inequality
$$
\sum_{\mu\in M_K} \log^-\max\{|r(u)|_\mu, |s(v)|_\mu\}
<-\epsilon\cdot (\max\{h(u),h(v)\}) \eqno(2.24)
$$
}
\medskip

\noindent {\it Proof }. We begin by observing that
the points $(u,v)$ such that $r(u)\cdot s(v)=0$ are contained in
a finite union of translates of proper subtori; hence we will tacitely
disregard these pairs. Also, for each constant $C$, the   pairs 
$(u,v)\in(\OS^\times)^2$ such that $\max\{h(u),h(v)\}\leq C$
 are finite in number,
so they certainly lie in finitely many translates of proper subtori.
For this reason, we can restrict our attention to ``large" solutions
of the inequality (2.23) or (2.24).
Let us write 
$$
\eqalign{
r(X)&=r_0X^a(X-\theta_1)^{a_1}\cdots(X-\theta_m)^{a_m},\cr
s(X)&=s_0X^b(X-\eta_1)^{b_1}\cdots(X-\eta_n)^{b_n};
}
$$
here $a$ and $b$ are non-negative integers with $ab=0$,
$a_1,\ldots,a_m,$ $b_1,\ldots,b_n$ are positive integers and
$\theta_1,\ldots,\theta_m,$ ({\it resp.} $\eta_1,\ldots,\eta_n$) are
pairwise distinct nonzero algebraic numbers in a number field $K$.
We shall first prove the first part of the Proposition,
so we shall estimate the relevant sum for $\mu$ running 
in the complement of $S$.   
Let, for each place $\mu\in M_K\setminus S$,
 $\delta_\mu$ be  the minimum of the set
$$
\eqalign{
&\{
|\theta_i-\theta_j|_\mu, {\rm \ for\ }
(i,j)\in\{1,\ldots,m\}\times\{1,\ldots,m\}, (i\neq j)\}\cr
&\cup \{|\eta_i-\eta_j|_\mu , {\rm \ for\ }
 (i,j)\in\{1,\ldots,n\}\times\{1,\ldots,n\}, (i\neq j)\} 
\cup \{|\theta_1|_\mu,\ldots,|\theta_m|_\mu,
|\eta_1|_\mu,\ldots,|\eta_n|_\mu\}.
}
$$
 Clearly, $\delta_\mu>0$ for all $\mu$ and $\delta_\mu=1$ for all
but finitely many places.  Let us divide the places 
of $M_K\setminus S$ in two classes ${\cal A}$ and ${\cal B}$:
 namely,
$$
{\cal A}:=\{\mu\in M_K\, |\, {\rm there\ exist\ } i\in\{1,\ldots, m\}
\ {\rm and}\ j\in\{1,\ldots,n\}\ {\rm such\ that }\ 
|u-\theta_i|_\mu<{\delta_\mu\over 2}\ {\rm and}\ |v-\eta_j|_\mu<{\delta_\mu\over 2}\}
$$
while ${\cal B}$ is its complement in the subset of $M_K\setminus S$
composed by the places $\mu$ such that $\max\{|r(u)|_\mu,|s(v)|_\mu\}<1$.
Then the left side of (2.23) becomes  
$$
\sum_{\mu\in M_K\setminus S}\log^-(\max\{|r(u)|_\mu, |s(v)|_\mu\})=
\sum_{\mu\in{\cal A}}\log^-(\max\{|r(u)|_\mu, |s(v)|_\mu\})+
\sum_{\mu\in{\cal B}}\log^-(\max\{|r(u)|_\mu, |s(v)|_\mu\}).
$$
The second term is simply bounded independently of $u,v$ by
$$
\sum_{\mu\in{\cal B}}\log^-(\max\{|r(u)|_\mu, |s(v)|_\mu\})\geq
-h(r_0)-h(s_0)+\max\{\deg r,\deg s\}\sum_{\mu\in\cal B}\log^{-}{\delta_\mu\over 2}.
\eqno(2.25)
$$
To estimate the first term, let us notice that for each $\mu\in{\cal A}$
there exists exactely one pair
$(i_\mu,j_\mu)\in\{1,\ldots,m\}\times\{1,\ldots,n\}$ such that
$|u-\theta_{i_\mu}|_\mu<{\delta_\mu\over 2}$, and $|v-\eta_{j_\mu}|_\mu<{\delta_\mu\over 2}$.
Then, for each such place $\mu$, we can bound
$$
\eqalign{
\log^{-}\max\{|r(u)|_\mu,|s(v)|_\mu\}&\geq
\log^{-}|r_0|_\mu+\log^{-}|s_0|_\mu
\cr &+\max\{\deg r,\deg s\}
\log^{-}{\delta_\mu\over 2}+\log^{-}
\max\{|u-\theta_{i_\mu}|_\mu^{a_{i_\mu}},|v-\eta_{j_\mu}|_\mu^{b_{j_\mu}}\},
}\eqno(2.26)
$$
where we have used the fact that $|u|^a_\mu=|v|^b_\mu=1$ for $\mu\in M_K\setminus S$.
Again the quantity   
$$
-\sum_{\mu\in {\cal A}}
(\log^{-}|r_0|_\mu+\log^{-}|s_0|_\mu+\max\{\deg r,\deg s\}
\log^{-}{\delta_\mu\over 2})
$$
is bounded independently of $u,v$ simply by $h(r_0)+h(s_0)+
\max\{\deg(r),\deg(s)\}
\sum_{\mu\in{\cal A}}\log^{-}{\delta_\mu\over 2}$.
It remains to estimate the term
$$
\sum_{\mu\in{\cal A}}\log^{-}\max\{|u-\theta_{i_\mu}|_\mu^{a_{i_\mu}},
|v-\eta_{j_\mu}|_\mu^{b_{j_\mu}}\}
$$
in the quantity appearing in (2.26). This can be rewritten as
$$
\sum_{(i,j)\in\{1,\ldots,m\}\times\{1,\ldots,n\}}\sum_{\mu\in {\cal A}_{i,j}}
\log^{-}\max\{
|u-\theta_{i}|^{a_{i_\mu}}_\mu,|v-\eta_{j}|^{b_{j_\mu}}_\mu\}
$$
where we define ${\cal A}_{i,j}$ to be the subset of ${\cal A}$
composed of the places $\mu$ for which $(i_\mu,j_\mu)=(i,j)$.
Proposition 3, applied with 
$\epsilon/(2(\deg r)(\deg s))$ instead of $\epsilon$, gives  for each pair $(i,j)$ 
the upper bound
$$
\sum_{\mu\in {\cal A}_{i,j}}
\log^{-}\max\{|u-\theta_{i}|_\mu,|v-\eta_{j}|_\mu\}
\geq -{\epsilon\over 2(\deg r)(\deg s)}\cdot \max\{h(u),h(v)\}\eqno(2.27)
$$
for $(u,v)\in(\OS^\times)^2$
outside a finite union of translates of proper subgroups, of which those of
dimension one can be effectively determined. Then, outside the union
of such exceptional subvarieties, we get from (2.25), (2.26) and (2.27),
$$
\sum_{\mu\in M_K\setminus S}\log^{-}\max\{|r(u)|_\mu,|s(v)|_\mu\}\geq
-{\epsilon\over 2}\max\{h(u),h(v)\}+O(1).\eqno(2.28)
$$
 This clearly implies the first part of Proposition 4.
We shall now estimate the relevant sum for $\mu$ running  over
the places in $S$. Recall that by hypothesis one at least between
the exponents $a$ and $b$ vanishes. Suppose that $b=0$, say.
We shall apply  Ridout's theorem (i.e. the case $N=2$ of the Subspace Theorem)
to give a lower bound for $\log^-|s(v)|_\mu$, valid for all but finitely many
$v$ (of course, in our situation, 
we could also apply the stronger Baker's bound for linear forms in logarithms).
For each $\mu\in S$, we denote by $j_\mu\in\{1,\ldots,n\}$ an index such that
for all $j=1,\ldots,n$
$$
|v-\eta_{j_\mu}|_\mu\leq |v-\eta_j|_\mu.
$$
Then
$$
|s(v)|_\mu\geq |s_0|_\mu\cdot |v-\eta_{j_\mu}|_\mu^{\deg s}.
$$
By Ridout's theorem we have that for all but finitely many $S$-units
$v\in \OS^\times$,
$$
|v-\eta_{j_\mu}|_\mu>H(v)^{-{\epsilon\over 3(\deg s)|S|}}\geq
\max\{H(u),H(v)\}^{-{\epsilon\over 3(\deg s)|S|}}
$$
so that, for all but finitely many $v\in\OS^\times$,
$$
\sum_{\mu\in S}\log^-\max\{|r(u)|_\mu,|s(v)|_\mu\}\geq
\sum_{\mu\in S}\log^- |s(v)|_\mu\geq -{\epsilon\over 3}\max\{h(u),h(v)\}+O(1).\eqno(2.29)
$$
Putting together (2.28) and (2.29) we obtain that outside the union of a finite set
and a finite number of effectively computable translates of one dimensional subgroups
we have
$$
\sum_{\mu\in S}\log^-\max\{|r(u)|_\mu,|s(v)|_\mu\}\geq -
{5\over 6}\epsilon\cdot  \max\{h(u),h(v)\}+O(1) 
$$
from which the Proposition follows immediately.

\bigskip

{\it Proof of Theorem 1}. We start by remarking that it suffices to prove  that the
relevant set is contained in a finite union of translates of proper subtori of
$\Gm$. This is very easy to see directly, since we are working in dimension 2; in any
case, this equivalence follows from a well-known theorem of M. Laurent 
(see e.g. [13]).
This remark applies also to the Main Theorem.\smallskip

To go on with the proof, we may clearly assume that the polynomials  
$p(X,Y), q(X,Y)$ are defined over a number field $K$ and that $\Gamma$ is the group
$(\OS^\times)^2$ for a suitable finite subset $S$ of $M_K$. 

We first slightly transform the relevant inequality. Observe that, by definition, 
$$
\eqalign{h(p(u,v):q(u,v):1)&=\sum_{\mu\in M_K}\log \max\{|p(u,v)|_\mu,
|q(u,v)|_\mu,1\}\cr &=\sum_{\mu\in M_K}\max\{0,\log \max\{ |p(u,v)|_\mu,
|q(u,v)|_\mu\}\}.}
$$
Then,  using the formula $\max\{0,\log(\cdot)\}=\log^+(\cdot)=
\log(\cdot)-\log^-(\cdot)$, we can write
$$
\eqalign{
h(p(u,v):q(u,v):1)&=\sum_{\mu\in M_K}\log \max\{|p(u,v)|_\mu, |q(u,v)|_\mu\}
-\sum_{\mu\in M_K}\log^-\max\{|p(u,v)|_\mu, |q(u,v)|_\mu\}\cr
&=h(p(u,v):q(u,v))-\sum_{\mu\in M_K}\log^-\max\{|p(u,v)|_\mu, |q(u,v)|_\mu\}.
}
$$
Note that $h(p(u,v)/q(u,v))=h(p(u,v):q(u,v))$ so the Theorem
is reduced to proving that the solutions to the inequality
$$
\sum_{\mu\in M_K}\log^-\max\{|p(u,v)|_\mu, |q(u,v)|_\mu\}<-\epsilon
\max\{h(u),h(v)\} \eqno(2.30)
$$
are contained in finitely many translates of proper subtori of $\Gm$,
of which those of dimension one can be effectively determined.
We shall prove this contention, by estimating separately the two subsums
$\sum_{\mu\in M_K\setminus S}\log^-\max\{|p(u,v)|_\mu, |q(u,v)|_\mu\}$
and $\quad \sum_{\mu\in S}\log^-\max\{|p(u,v)|_\mu, |q(u,v)|_\mu\}$.

Let us begin with the second.  

By hypothesis not both $p(X,Y)$
and $q(X,Y)$ vanish at the origin; suppose, for instance,
that $p(0,0)\neq 0$. 
 We apply  Proposition 1, with $f(X,Y)=p(X,Y)$  
and ${\epsilon/2}$ instead of $\epsilon$.
Then (the second part of) Proposition 1 states that, outside the union of a finite
set and finitely many effectively computable translates of one 
dimensional subgroups, one has  
$$
\sum_{\mu\in S}\log^- |p(u,v)|_\mu>
-{\epsilon\over 2}\max\{h(u),h(v)\} 
$$
so clearly  
$$
\sum_{\mu\in S}\log^-\max\{|p(u,v)|_\mu, |q(u,v)|_\mu\}>-
{\epsilon\over 2}\max\{h(u),h(v)\}. \eqno(2.31)
$$

To estimate the sum for $\mu$ running over the complement of $S$, we argue as follows:
Since $p(X,Y), q(X,Y)$ are co\-prime polynomials, there exist
polyno\-mials $A(X,Y)$, $B(X,Y),C(X,Y)$, $D(X,Y)\in\OS[X,Y]$ 
and nonzero polynomials $r(X),s(Y)$ such that 
$$
\eqalign{
A(X,Y)p(X,Y)+B(X,Y)q(X,Y)&=r(X)\cr
C(X,Y)p(X,Y)+D(X,Y)q(X,Y)&=s(Y)
}
$$
Clearly, for a place $\mu\not\in S$ and $S$-unit point $(u,v)$,
we can bound by $1$ the quantities $|A(u,v)|_\mu,|B(u,v)|_\mu$,
$|C(u,v)|_\mu, |D(u,v)|_\mu$.
Then
$$
\sum_{\mu\in M_K\setminus S}\log^-\max\{|p(u,v)|_\mu,|q(u,v)|_\mu\}\geq 
\sum_{\mu\in M_K\setminus S}
\log^-\max\{|r(u)|_\mu,|s(v)|_\mu\}.
$$
The last quantity can be bounded, using Proposition 4,
by $-(\epsilon/2)\max\{h(u),h(v)\}$ on the 
$S$-unit points $(u,v)$ outside a finite union of  
translates of proper subgroups of $\Gm$. This inequality, combined with (2.31),
concludes the proof of Theorem 1.

\bigskip

\noindent{\it Proof of Corollary 1}. We could derive Corollary 1 both from Theorem 1 and from
Proposition 2. We choose the latter possibility. We always have 
$$
\eqalign{
h((u-1)/(v-1))&=h(u-1:v-1)=\sum_{\mu\in M_K}\log\max\{|u-1|_v,|v-1|_\mu\}\cr
&\leq \sum_{\mu\in M_K}\log\max\{1,|u-1|_v,|v-1|_\mu\}=h(1:u-1:v-1).
}
$$
Also, $h(1:u-1:v-1)\leq h(1:u:v)+O(1)$. So in particular
$$
\lim\sup_{(u,v)} {h((u-1)/(v-1))\over h(1:u:v)}\leq 1.
$$
We shall prove the non trivial inequality 
$$
\lim\inf_{(u,v)\in\Gamma} {h((u-1/(v-1))\over h(1:u:v)}\geq 1,
$$
where now the limit is taken only over the multiplicative independent 
pairs $(u,v)\in\Gamma$. This means that given a positive $\epsilon$, only 
finitely many such pairs should verify the inequality
$$
h((u-1)/(v-1))<(1-\epsilon)h(1:u:v).
$$
As we have already remarked,
$$
h((u-1)/(v-1))\geq h(1:u-1:v-1)+\sum_{\mu\in M_K}\log^-\max\{|u-1|_\mu,|v-1|_\mu\}
$$
and by Proposition 2 all but finitely many
multiplicative independent pairs $(u,v)\in\Gamma$ satisfy
$$
\sum_{\mu\in M_K}\log^-\max\{|u-1|_\mu,|v-1|_\mu\}
\geq -{\epsilon}\max\{h(u),h(v)\}\geq
-{\epsilon}h(1:u:v).
$$

\bigskip

\noindent{\it Proof of Main Theorem}. Let us write 
$f(X,Y)=p(X,Y)/q(X,Y)$ for polynomials $p(X,Y),q(X,Y)$
defined over a number field $K$. Note that the hypothesis
that $1\in\{T_1(X,Y),\ldots,T_N(X,Y)\}$ appearing on the 
Main Theorem is equivalent to the hypothesis in Theorem 1
that $p(X,Y),q(X,Y)$ do not both vanish at $(0,0)$.  
 The Main Theorem then
follows essentially from Proposition 1 and  Theorem 1.
We have just to compare the maximum of the heights of the
monomials appearing  in the development of $p(X,Y),q(X,Y)$
with the maximum of the heights of $u$ and $v$.
This is the object of Lemma 2.
 
First of all, note that (the elementary) Lemma 1 immediately implies the
result in a stronger form if all the non constant  monomials $T_i(X,Y)$ 
vanish on a same one-dimensional subtorus. In fact, in such a case, the
function $f(X,Y)$ could be written as $f(X,Y)=\varphi(X^aY^b)$ for 
a suitable one-variable rational function $\varphi(T)\in K(T)$ and a
pair $(a,b)\in\Z^2$; in that case the proof is exactely the same as in the
corresponding  case of the proof of Proposition 1. 
 So we suppose that the set of monomials 
$\{T_1,\ldots,T_N\}\setminus\{1\}$ satisfies the hypothesis of Lemma 2.

Fix a positive $\epsilon$.
Since $h(f(u,v))=h(p(u,v):q(u,v))$, we have
$$
h(f(u,v))\geq \max\{h(p(u,v)),h(q(u,v))\}+\sum_{\mu\in M_K}\log^{-} \max\{|p(u,v)|_\mu,
|q(u,v)|_\mu\};\eqno(2.32)
$$
by Proposition 1 the first summand is bounded from below as
$$
\max\{h(p(u,v)),h(q(u,v))\}\geq
(1-\epsilon/ 2)\max\{h(T_1(u,v)),\ldots,h(T_N(u,v))\}\eqno(2.33)
$$ 
outside a finite union of translates of proper subgroups. 
From the proof of Theorem 1 ({\it see} inequality (2.31)),
the second summand in (2.32) is bounded, outside a finite union of
translates of proper subgroups of $\Gm$, as
$$
 \sum_{\mu\in M_K}\log^{-}\max\{|p(u,v)|_\mu, |q(u,v)|_\mu\}>-
{\epsilon\over 4\max\{\deg_X p,\deg_Y p,\deg_X q,\deg_Y q\}}
\max\{h(u),h(v)\}
$$
so by Lemma 2 we have  
$$
\sum_{\mu\in M_K}\log^{-}\max\{|p(u,v)|_\mu, |q(u,v)|_\mu\}
>-{\epsilon\over 2}\max\{h(T_1(u,v)),\ldots,h(T_N(u,v))\}.
\eqno(2.34)
$$
Putting together (2.32), (2.33) and (2.34) we obtain that
the solutions to (1.4) are contained in a finite union of proper
translates, as wanted.
Using  this fact and Lemma 2 we obtain (1.5).

\bigskip

{\bf References}
\bigskip

\item{[1]} Y. Bugeaud, P. Corvaja, U. Zannier, An upper bound for the
G.C.D. of $a^n-1$ and $b^n-1$, {\it Math. Zeit.} {\bf 243} (2003), 79-84.
\smallskip

\item{[2]} Y. Bugeaud, F. Luca, A quantitative upper bound for the greatest
prime factor of $(ab+1)(ac+1)(bc+1)$, to appear in {\it Acta Arith.}.
\smallskip

\item{[3]} C. Fuchs, An upper bound for the G.C.D. of two linear recurring
 sequences, {\it Math. Slovaca} {\bf 53} (2003), 21-42.\smallskip

\item{[4]} S. Hernandez, F. Luca, On the largest prime factor of 
$(ab+1)(ac+1)(bc+1)$, to appear in {\it Bol. Soc. Mat. Mexicana}.
\smallskip

\item{[5]} P. Corvaja, U. Zannier,  Diophantine equations with power sums and universal
Hilbert sets, {\it Indagationes Math.}, N.S. {\bf 9 (3)} (1998), 317-332.
\smallskip

\item{[6]}  P. Corvaja, U. Zannier, 
On the greatest prime factor of $(ab+1)(ac+1)$.  
{\it Proceedings American Math. Soc.} {\bf 131,6} (2003),
1705-1709.
\smallskip

\item{[7]} P. Corvaja, U. Zannier, Finiteness of integral values for the ratio of two linear
recurrences, {\it Inv. Math} {\bf 149} (2002), 431-451.\smallskip

\item{[8]}  J.H. Evertse, On sums of $S$-units and linear recurrences,
{\it Compositio Math.}, {\bf 53} (1984), 225-244.
\smallskip

\item{[9]}  M. Hindry, J.H. Silverman, {\it Diophantine Geometry},
Springer-Verlag, 2000.\smallskip

\item{[10]} S. Lang  Fundamentals of Diophantine Geometry, {\it Graduate
Texts in Mathematics} {\bf }, Springer Verlag.
\smallskip

\item{[11]}  W.M. Schmidt, {\it Diophantine Approximation}, Springer-Verlag LNM
785.\smallskip

\item{[12]}  W.M. Schmidt, {\it Diophantine Approximations and Diophantine
Equations}, Springer-Verlag LNM 1467, 1991.\smallskip 

\item{[13]}  U. Zannier, {\it Some applications of Diophantine
Approximation to Diophantine Equations}, Forum Editrice, Udine (2003).\smallskip

\bigskip

\vfill

P. Corvaja\hfill U. Zannier

Dip. di Matematica e Informatica\hfill  Ist. Univ. Arch. - D.C.A.

Via delle Scienze\hfill S. Croce, 191

33100 - Udine (ITALY)\hfill  30135 Venezia (ITALY)

corvaja@dimi.uniud.it\hfill zannier@iuav.it

\bye